%
%
\input amstex.tex
\documentstyle{amsppt}
\hsize=5.3in
\parindent=2pc
\TagsOnRight
%
%
\def\nologo{\expandafter\let\csname logo\string @\endcsname=\empty}
\def\section#1{\par\bigpagebreak
    \csname subhead\endcsname #1\endsubhead\par\medpagebreak}
\def\subsection#1{\par\medpagebreak
    \csname subsubhead\endcsname #1 \endsubsubhead}
\def\Theorem#1#2{\csname proclaim\endcsname{Theorem #1} #2 
    \endproclaim}
\def\Corollary#1{\csname proclaim\endcsname{Corollary} #1 \endproclaim}
\def\Proposition#1#2{\csname proclaim\endcsname{Proposition #1} #2 
    \endproclaim}
\def\Lemma#1#2{\csname proclaim\endcsname{Lemma #1} #2 \endproclaim}
\def\Remark#1#2{\remark{Remark {\rm #1}} #2 \endremark}

\def\C{{\Bbb C\,}} 

\def\R{{\Bbb R}}
\def\Dt{\Delta}
\def\dt{\delta}
\def\vep{\varepsilon}
\def\ld{\lambda}
\def\ox{\otimes}
\def\oxC{\otimes_\C}
\def\half{\frac{1}{2}}
\def\U{{\Cal U}}
\def\V{{\Cal V}}
\def\A{{\Cal A}}

\def\H{{\Cal H}}
\def\Uq{U_q(\frak{g})}
\def\Uqk{U_q^{\text{tw}}(\frak{k})}
\def\Aq{A_q(G)}
\def\Aq{A_q(G)}
\def\Aq{A_q(G)}

\def\HomC{\operatorname{Hom}_\C}

\def\EndC{\operatorname{End}_\C}
\def\id{\operatorname{id}}

\NoRunningHeads
\document
\par\noindent
{\eightpoint
{\it To appear in}: Proceedings of the XXth International 
Colloquium on Group Theoretical Methods in Physics (ICGTMP),
July 4--9, Toyonaka, Japan.
}
\par\bigskip\par\bigskip
\centerline{\bf QUANTUM SYMMETRIC SPACES AND}
\centerline{\bf RELATED $q$-ORTHOGONAL POLYNOMIALS}
\par\bigskip
\par\bigskip
\centerline{MASATOSHI NOUMI}
\par\smallskip
\centerline{and}
\par\smallskip
\centerline{TETSUYA SUGITANI}
\par\smallskip
{\eightpoint\it
\centerline{Department of Mathematical Sciences, University of Tokyo}
\centerline{Komaba 3-8-1, Meguro-Ku, Tokyo 153, Japan}
} 
\par\bigskip
\par\bigskip
{\eightpoint\rm 
\centerline{ABSTRACT}
\midinsert
\narrower\narrower
\noindent
A class of quantum analogues of compact symmetric spaces of 
classical type is introduced 
by means of constant solutions to the reflection 
equations.  Their zonal spherical functions are discussed 
in connection with $q$-orthogonal polynomials. 
\endinsert
}
\par\bigskip
\tenpoint
\noindent
The following two naive questions are the main motives of this
paper: 
\roster
\item How can one define the analogue of homogeneous spaces $G/K$ 
in the framework of quantum groups? 
\item What do their spherical functions look like?
\endroster
To be more specific, we will introduce a class of 
quantum analogues of compact symmetric spaces $G/K$ of classical
type and study their 
zonal spherical functions associated with finite 
dimensional representations.  
It is natural to expect that they could provide a good class of 
$q$-orthogonal polynomials in many variables. 
Unfortunately, we have not yet reached an abstract definition of 
quantum symmetric spaces.  In this paper we will propose
instead a practical method to construct examples of quantum 
symmetric spaces of classical type, by means of constant 
solutions to the reflection equations. 
This method works well in fact and it turns out in many 
examples that the zonal spherical functions are expressed by 
the Macdonald polynomials associated with root systems or 
by Koornwinder's Askey-Wilson polynomials for $BC_\ell$.
\section{1. Recalls on compact quantum groups}
Let $G$ be one of the compact classical groups $SU(N)$, $SO(N)$, 
$Sp(N)$ and let $\frak{g}$ be the complexification 
of the Lie algebra of $G$. 
For such a compact group $G$, we already have a (more or less) 
standard definition of the quantum group $G_q$
(Woronowicz${}^1$, Reshetikhin-Takhtajan-Faddeev${}^2$,
see also Hayashi${^3}$, Dijkhuizen-Koornwinder${}^4$).  
One can define in fact two Hopf $\ast$-algebras (over $\C$) 
$\Uq$ and $\Aq$ which are $q$-deformations of 
the universal enveloping algebra $U(\frak{g})$ of $\frak{g}$ and 
of the algebra of regular functions $A(G)$ on $G$, respectively.
For $\Uq$, we take the quantized universal enveloping 
algebra of Drinfeld and Jimbo endowed with the $\ast$-operation
corresponding to the compact real form. 
Hereafter we always assume that 
$q$ is a real number with $0<q<1$.  The algebra 
$\Aq$ is the subalgebra of the algebraic dual 
$\Uq^{\vee}=\HomC(\Uq,\C)$ generated by the 
matrix elements of the vector representation.  
Hence we have the natural pairing of Hopf $\ast$-algebras
$$
(\ ,\ ) : \Uq\times\Aq \to \C. \tag{1.1}
$$
\par
We say that a Hopf algebra $\A$ is a {\it Hopf $\ast$-algebra} 
if it has a $\ast$-operation (involutive, conjugate linear 
antiautomorphism of $\R$-algebra) such that 
the coproduct $\Dt:\A\oxC\A\to\A$ and the counit $\vep: \A\to\C$ 
are $\ast$-homomorphisms. 
When $\U$ and $\A$ are two Hopf $\ast$-algebras, 
a $\C$-bilinear mapping $(\ ,\ ): \U\times\A\to\C$ is called 
a {\it pairing of Hopf $\ast$-algebras} if the following 
conditions are satisfied:
\roster
\item $(ab,\varphi)=(a\ox b, \Dt(\varphi)),
\quad (1,\varphi)=\vep(\varphi)$, 
\item $(a,\varphi\psi)=(\Dt(a),\varphi\ox\psi),
\quad (a,1)=\vep(a)$, 
\item $(S(a),\varphi)=(a,S(\varphi))$,
\item $(a^\ast,\varphi)=\overline{(a,S(\varphi)^\ast)}$,
\endroster
for any $a,b\in\U$ and $\varphi,\psi\in\A$. 
\par
{}From the pairing between $\Uq$ and $\Aq$ we obtain 
a natural structure of $\Uq$-bimodule on $\Aq$, 
corresponding to the right and the left regular representations
of $G$.  
The left and the right actions of an element $a\in\Uq$ on $\Aq$ are 
defined by
$$
a.\varphi=(a,\Dt(\varphi))_2,\quad 
\varphi.a=(a,\Dt(\varphi))_1 
\quad(\varphi\in\Aq),\tag{1.2}
$$
where $(a,\,\cdot)_i$ stands for the contraction with respect 
to the $i$-th tensor component. 
It should be noted here that these actions of $\Uq$ on $\Aq$ 
are compatible with the multiplication of $\Aq$: If $a\in\Uq$ 
and $\Dt(a)=\sum_{(a)} a_{(1)}\ox a_{(2)}$, then we have
$$
a.(\varphi\psi)=\sum_{(a)}(a_{(1)}.\varphi) \,\,(a_{(2)}.\psi),\quad
(\varphi\psi).a=\sum_{(a)}(\varphi.a_{(1)})\,\,(\psi.a_{(2)}), \tag{1.3}
$$
for any $\varphi,\psi\in\Aq$ and $a.1=1.a=\vep(a)1$. 
\par
The most fundamental fact about the structure of 
the $\Uq$-bimodule $\Aq$ is that $\Aq$ has the 
the following irreducible decomposition of Peter-Weyl type:
$$
\Aq=\bigoplus_{\ld\in P^+_G} W(\ld),
\quad W(\ld)\simeq V(\ld)^\vee\oxC V(\ld). \tag{1.4}
$$
Here $P^+_G$ denotes the cone of dominant integral weights 
corresponding to the $G$-rational representations. 
For each $\ld\in P^+_G$, $V(\ld)$ is the finite dimensional 
irreducible left $\Uq$-module (or right $\Aq$-comodule) 
with highest weight $\ld$, and $W(\ld)$ is the vector 
subspace of $\Aq$ spanned by the matrix elements of $V(\ld)$.
\par
Another characterization of $W(\ld)$ is given by the 
action of the center ${\Cal Z}\Uq$ of $\Uq$.  For each 
$\ld\in P^+_G$, we denote the central character of $V(\ld)$
by $\chi_\ld:{\Cal Z}\Uq\to\C$: 
$C|_{V(\ld)}=\chi_\ld(C)\id_{V(\ld)}\ \ (C\in{\Cal Z}\Uq)$. 
Then $W(\ld)$ is the following simultaneous eigenspace of 
${\Cal Z}\Uq$:
$$
W(\ld)=\{ \varphi\in\Aq\,|\, C.\varphi=\chi_\ld(C)\varphi 
\ \text{for all}\ C\in {\Cal Z}\Uq\}.\tag{1.5}
$$
\par
We denote by $h:\Aq\to W(0)=\C$ the projection to the 
trivial representation in the decomposition (1.4).  Then 
the functional $h$ gives the unique invariant functional
with $h(1)=1$, which we call the {\it normalized Haar functional} 
of the quantum group $G_q$. 
By this functional, we can define a scalar product on $\Aq$:
$$
\langle{\varphi|\psi}\rangle=h(\varphi^\ast\,\psi)\quad
(\varphi,\psi\in\Aq). \tag{1.6}
$$
It is known that $\langle{\ | \  }\rangle$ is in fact 
a positive definite Hermitian form and that the Peter-Weyl 
decomposition (1.4) is orthogonal under this scalar product. 
\section{2. Quantum analogue of $G/K$ and $K\backslash G/K$}
In this section, we discuss how one can define the quantum analogue 
of homogeneous spaces $G/K$ and double coset spaces 
$K\backslash G/K$ in the sense of $q$-deformation of algebras 
of regular functions. 
Suppose now a closed subgroup $K$ of $G$ is given and let 
$\frak{k}$ be the complexification of the Lie algebra of $K$. 
The question we have to discuss first is: 
\roster
\item"$\bullet$"{\it What is the appropriate 
definition of a quantum subgroup?}
\endroster
We will take here the infinitesimal approach of quantized universal 
enveloping algebras rather than the global approach of quantized 
algebras of functions.  Given a pair $(\frak{g},\frak{k})$
of a Lie algebra $\frak{g}$ and its Lie subalgebra $\frak{k}$, 
it is natural to ask: 
\roster
\item"$\bullet$"
{\it What should be the quantum analogue 
of the pair $(U(\frak{g}),U(\frak{k}))$? }
\endroster 
We have at least two possible versions of quantum analogue of 
the pair $(U(\frak{g}),U(\frak{k}))$.  
\par\smallskip\noindent
(A) {\it Regular version\/} 
(or a quantum subgroup in the {\it strict} sense).
\newline
Suppose we have a Hopf subalgebra $\V$ of $\U=\Uq$ such that 
$\V$ ``tends" to $U(\frak{k})$ as $q\to 1$. In such a case
we might write $\V=U_q(\frak{k})$ as well and might  
expect also that there would be a quotient Hopf algebra, 
say $A_q(K)$, of $\Aq$ representing the quantum subgroup 
$K_q$. 
As for quantized universal enveloping algebras of Drinfeld-Jimbo, 
a class of pairs $(U_q(\frak{g}),U_q(\frak{k}))$ of Hopf algebras 
arises naturally from the embedding of root systems.  
This class is, however, not so large as to cover all symmetric pairs. 
In fact, no natural embedding $U_q(\frak{so}(N))\to U_q(\frak{sl}(N))$,
nor $U_q(\frak{sp}(N))\to U_q(\frak{sl}(N))$ seems to be known.
In relation to this point, non-existence of Hopf algebra 
homomorphisms $A_q(G)\to A_q(K)$ for pairs $(G,K)$ of 
classical groups is discussed by Hayashi${}^5$.
In our context, we should probably say that the condition
$$
\Dt(\V)\subset \V\ox\V \tag{2.1}
$$
for a Hopf subalgebra $\V\subset\U=\Uq$ is too restritive. 
\par\smallskip\noindent
(B) {\it Twisted version\/} (or a quantum subgroup in the 
{\it broader} sense).\newline 
We propose to consider subalgebras $\V\subset\U=\Uq$ which are 
not Hopf subalgebras, as well. Instead of (2.1), we assume that 
the subalgebra $\V$ satisfies the {\it coideal property}
$$
\Dt(\V)\subset \V\ox\U + \U\ox\V.   \tag{2.2} 
$$
This condition (2.2) is equivalent to saying that $\V$ is 
generated by a coideal of $\U$. 
One might call a subalgebra $\V$ of $\U$ satisfying (2.2) a 
{\it coideal subalgebra}. Note that, in this context of subalgebras, 
it is not natural to impose the condition on the counit. 
Subalgebras of this type are already considered by 
Olshanski${}^6$ in the case of Yangians under the name of  
{\it twisted Yangians}.  As for $\Uq$, 
the $q$-deformation of $\frak{so}(N)$ proposed by 
Gavrilik and Klimyk${}^7$ can be regarded as a 
subalgebra of $U_q(\frak{sl}(N))$ satisfying (2.2) 
(see Section 2.4 of Noumi${}^8$). 
We remark that, in these examples, 
the subalgebras in question are actually {\it one-sided} 
in the sense 
$$
\Dt(\V)\subset \V\ox\U\quad\text{or}\quad
\Dt(\V)\subset \U\ox\V.\tag{2.3}
$$
When we can consider a subalgebra $\V\subset\U=\Uq$ satisfying 
(2.2) as a $q$-deformation of $U(\frak{k})$, we will use the 
notation  $\V=U_q^{\text{tw}}(\frak{k})$ to remember that $\V$ 
is may not be a Hopf subalgebra. 
\par\smallskip\par
We now suppose that a subalgebra $\V=\Uqk$ 
of $\U=\Uq$ satisfying (2.2) is given and define the left 
ideal ${\Cal J}$ of $\U$ by 
$$
{\Cal J}=\sum_{a\in\V} \U(a-\vep(a)). \tag{2.4}
$$
Then 
we define the quantum analogue of the algebra of right 
$K$-invariant regular functions on $G$ to be 
the subspace of all elements in $\Aq$ annihilated 
by the left action of ${\Cal J}$. Namely, we set
$$
A_q(G/K):=\{\varphi\in\Aq\,|\, {\Cal J}.\varphi=0\}. \tag{2.5}
$$
The coideal property (2.2) of $\V=\Uqk$ garantees that this subspace
$A_q(G/K)$ is actually a subalgebra of $\Aq$. 
Note also that $A_q(G/K)$ is a right $\Uq$-submodule of $\Aq$. 
Assume furthermore that 
the left ideal ${\Cal J}$ associated with $\V$ 
satisfies the condition
$$
S({\Cal J})^{\ast}\subset {\Cal J}. \tag{2.6}
$$
Then one can show that the subalgebra $A_q(G/K)$ of (2.5) becomes 
a $\ast$-subalgebra of $\Aq$.
\par
We will mainly consider the case when the pair $(\U,\V)$ is a 
{\it Gelfand pair} in the sense that 
$$
\dim_\C V(\ld)^{\Cal J} \le 1
\quad\text{for all }\quad \ld\in P^+, \tag{2.7}
$$
where $V(\ld)^{\Cal J}=\{v\in V(\ld)\,|\,{\Cal J}.v=0\}$.
Then the Peter-Weyl decomposition (1.4) implies the 
following multiplicity free decomposition of $A_q(G/K)$ 
as a right $\Uq$-module:
$$
A_q(G/K)\simeq\bigoplus_{\ld\in P^+_{G,{\Cal J}}} V(\ld)^\vee,\tag{2.8}
$$
where $P^+_{G,{\Cal J}}$ is the subset of $P^+_G$ consiting 
of all $\ld$ for which $V(\ld)$ has nonzero ${\Cal J}$-fixed vectors. 
\par
The next object we have to consider is the quantum analogue of 
the double coset space
$K\backslash G/K$. 
By using the group-like element $q^\rho$ of $\Uq$, 
corresponding the half sum of positive roots, we modify  
the $\ast$-operation of $\Uq$ as 
$\overline{a}=q^{\rho} a^\ast q^{-\rho}.$
With this notation, we take the left ideal ${\Cal J}$ 
of (2.4) and the right ideal $\overline{\Cal J}$ of $\U=\Uq$ 
and set 
$$
A_q(K\backslash G/K):=
\{\varphi\in\Aq\,|\, {\Cal J}.\varphi=
\varphi.\overline{\Cal J}=0
\}. \tag{2.9}
$$
Under the condition (2.6), $A_q(K\backslash G/K)$ becomes 
a $\ast$-subalgebra of $\Aq$, which we regard as the quantum 
analogue of the algebra of $K$-biinvariant regular functions on $G$. 
In what follows, we also use the notation $\H=A_q(K\backslash G/K)$ 
for simplicity. Under the condition (2.7) of Gelfand pair, 
this algebra $\H=A_q(K\backslash G/K)$ is decomposed into 
{\it one-dimensional} subspaces as 
$$
\H=\bigoplus_{\ld\in P^+_{G,{\Cal J}}}\H(\ld),
\quad\H(\ld)=\H\cap W(\ld). \tag{2.10}
$$
{}From (1.5), we see that (2.10) also gives the decomposition of 
$\H$ into simultaneous eigenspaces of the action of 
the center ${\Cal Z}\Uq$.  
We call a nonzero element $\varphi_\ld$ of $\H(\ld)$ 
$(\ld\in P^+_{G,{\Cal J}})$ a {\it zonal spherical function
on the quantum homogeneous space $(G/K)_q$}, 
associated with the representation $V(\ld)$. In other words,
the zonal spherical function is characterized by the conditions
${\Cal J}.\varphi_\ld=\varphi_\ld.\overline{\Cal J}=0$ and
$$
C.\varphi_\ld=\chi_\ld(C)\varphi_\ld
\quad\text{for all}\quad C\in{\Cal Z}\Uq,
\tag{2.11}
$$
up to a constant multiple.
Note also that we have Schur's orthogonality 
relations
$$
\langle{\varphi_\ld | \varphi_\mu}\rangle = 0 \quad(\ld\ne\mu)\tag{2.12}
$$
for $\ld,\mu\in P^+_{G,{\Cal J}}$, under the scalar 
product (1.6) defined by the Haar functional. 
\par
At this stage, our problems can be stated as follows:
\roster
\item Find a systematic way to construct a coideal subalgebra 
$\Uqk$ of $\Uq$ for a given Gelfand pair $(\frak{g},\frak{k})$. 
\item Describe the structure of subalgebras  
$A_q(G/K)$ and $A_q(K\backslash G/K)$ of invariants in 
$\Aq$. 
\endroster 
In many examples, the subalgebra $\H=A_q(K\backslash G/K)$ 
turns out to be commutative.  
If we have a subalgebra $\Uqk$ of $\Uq$ 
with the desired properties and 
if we can describe the structure of $\H$, then the zonal spherical 
functions $\varphi_\ld$ ($\ld\in P^+_{G,{\Cal J}}$)
would give rise to a family of $q$-orthogonal 
polynomials in many variables. 
Schur's orthogonality relations (2.12) 
would then guarantee their orthogonality relations. 
Furthermore, the action of the center ${\Cal Z}\Uq$ on 
$\H$ would provide a commuting family of $q$-difference operators 
for which the $q$-orthogonal polynomials should be simultaneous 
eigenfunctions.  
We will see below that this is the case in many examples
and that zonal spherical functions on quantum homogeneous 
spaces generate in fact various $q$-orthogonal polynomials. 
\par
There have been a lot of works related to the quantum analogue
of symmetric spaces of rank one and their spherical functions. 
See Vaksman-Soibelman${}^9$, Masuda et al.${}^{10}$, 
Koornwinder ${}^{11}$, and also Vilenkin-Klimyk${}^{12}$, 
for the interpretation of little $q$-Jacobi polynomials as 
the the matrix elements of unitary representations of $SU_q(2)$. 
Several extensions of this result are discussed by 
Noumi-Mimachi${}^{13}$, Koornwinder${}^{14}$, 
Koelink${}^{15}$ in relation to Podles's quantum spheres 
and Askey-Wilson polynomials. 
For spherical functions related to $SU_q(2)$, we 
refer to the survey papers Koornwinder${}^{16}$ and 
Noumi${}^{17}$.
For spherical functions on higher dimensional quantum spheres,
see Noumi-Yamada-Mimachi${}^{18}$, Vaksman-Soibelman${}^{19}$
and Sugitani${}^{20}$. 
\par
In the rank one case, a typical example of the 
{\it regular version} is the quantum sphere 
$(SU(2)/S(U(1)\times U(1)))_q$ 
obtained as the quotient space of $SU_q(2)$ by 
its diagonal subgroup, 
whose zonal spherical functions are expressed by little 
$q$-Legendre polynomials.  The first example of the 
{\it twisted version} 
is $(SU(2)/SO(2))_q$ and its zonal spherical functions are 
expressed then by continuous $q$-Legendre polynomials.   
It should be noted that both 
$(SU(2)/S(U(1)\times U(1)))_q$ and $(SU(2)/SO(2))_q$ are 
two extreme cases of quantum 2-spheres of Podles${}^{21}$. 
The zonal spherical functions on $(SU(2)/SO(2))_q$ were 
discussed by Koornwinder${}^{14}$ for the first time, 
in which he used the infinitesimal approach of twisted primitive 
elements in $U_q(\frak{sl}(2))$.  
This work of Koornwinder can be regarded as the starting point 
of quantum subgroups in the broader sense that we described above. 
\par
{}From the examples studied so far, it seems natural to suspect that 
quantum homogeneous spaces of the regular version would  
provide $q$-orthogonal polynomials with respect to discrete 
measures of $q$-integral.  On the other hand, quantum homogeneous 
spaces of the twisted version would be related in general 
to $q$-orthogonal polynomials with respect to mesures involving  
continuous parts. 
The quantum analogues of symmetric spaces that we discuss below 
belong to this second category. 
\section{3. Quantum symmetric spaces and reflection equations}
Recall first we have the following seven series of compact 
Riemannian symmetric spaces $G/K$ of classical type
(see Loos${}^{22}$): 
$$
\align
\text{AI:} & \quad SU(n)/SO(n)  \quad(N=n), \cr
\text{AII:} & \quad SU(2n)/Sp(2n) \quad(N=2n), \cr
\text{AIII:} & \quad U(n)/U(\ell)\times U(n-\ell)
\quad(N=n,\ell\le[\frac{n}{2}]),\cr
\text{BDI:} & \quad SO(N)/SO(\ell)\times SO(N-\ell)
  \quad(N=2n+1\ \ \text{or}\ \ 2n,\ell\le[\frac{N}{2}]),\tag{3.1}\cr
\text{CI:} & \quad Sp(2n)/U(n)\quad(N=2n),\cr
\text{CII:} & \quad Sp(2n)/Sp(2\ell)\times Sp(2(n-\ell))
\quad(N=2n,\ell\le[\frac{n}{2}]),\cr
\text{DIII:} & \quad SO(2n)/U(n) \quad(N=2n).
\endalign
$$ 
In this section, we will propose a method to construct 
{\it twisted} quantized universal enveloping algebra $\Uqk$
for $G/K$ of these series of symmetric spaces other than AIII.
The procedure of our ``quantization'' can be 
described in a unified manner by using constant solutions 
of the reflection equations. 
Although we will not treat the case of 
type AIII in this paper, a similar argument can be carried 
out as well, by using a different type of reflection equations.
We will not consider here the cases of group manifolds $G\times G/\Dt(G)$ 
either. 
\par
We remark that the quantization of symmetric spaces of type AI 
was discussed by Ueno-Takebayashi${}^{23}$ and, on the quantum 
$SU(3)/SO(3)$, they gave the interpretation of Macdonald 
polynomials of type $A_2$ as zonal spherical functions. 
In our previous paper ${}^8$ we discussed the quantization 
of AI and AII as well as the relation with 
the Macdonald polynomials of type $A$.  
The method we use here is essentially the same as that 
of the paper${}^8$
\par
Let $(G,K)$ be the pair of compact Lie groups corresponding 
to one of the symmetric space $G/K$ listed in (3.1).    
We assume that $G/K$ is not of type AIII. 
We denote by $V=\C^N$ the underlying vector space of the 
vector representation of $G$. 
We will use the following $R$-matrix $R\in\EndC(V\oxC V)$ 
for the vector representation of $\Uq$:
$$
R=q^{-1/N}(\sum_{1\le i,j\le N}e_{ii}\ox e_{jj}q^{\dt_{ij}}
+(q-q^{-1})\sum_{1\le j<i\le N} e_{ij}\ox e_{ji})\tag{3.2}
$$
if $G=SU(N)$, and
$$
R=\sum_{1\le i,j\le N}e_{ii}\ox e_{jj}q^{\dt_{ij}-\dt_{ij'}}
+(q-q^{-1})\sum_{1\le j<i\le N}
(e_{ij}\ox e_{ji} -\kappa_i\kappa_j q^{\rho_{i}-\rho_j}e_{ij}\ox e_{i'j'})
\tag{3.3}
$$
if $G=SO(N)$ or $G=Sp(N)$, where the $e_{ij}$'s are the matrix units 
corresponding to the canonical basis of $V$. 
For $G=SO(N), Sp(N)$, we use the notation 
$j'=N+1-j$ and $q^{\rho_j}$ is the $(j,j)$-component of the 
diagonal matrix representing $q^\rho\in\Uq$ 
on the vector representation $V$. 
If $G=SO(N)$, then $\kappa_j=1$ for all $j$, and 
if $G=Sp(N)$ ($N=2n$), then $\kappa_j=1$ or $-1$ according 
as $j\le n$ or $j>n$.  
We also use the $L$-operators ${\Cal L}^+,{\Cal L}^-\in\EndC(V)\ox\Uq$ 
such that 
$$
(\id_V\ox\rho_V)({\Cal L}^{\pm})=R^{\pm}\quad 
\text{with}\quad
R^+=PRP, \ R^{-}=R^{-1},
\tag{3.4}
$$
where $P\in\EndC(V\oxC V)$ is the flip $u\ox v\mapsto v\ox u$.
\par
The first step of our method is to find an appropriate 
constant solution to the following reflection equation 
for $J\in\EndC(V)$:
$$
R_{12}\,J_1 \,R_{12}^{t_1}\,J_2=J_2\,R_{12}^{t_1}\,J_1\,R_{12},
\tag{3.5}
$$
where $R_{12}^{t_1}$ stands for the matrix obtained from 
$R=R_{12}$ by transposition in the first tensor component. 
Reflection equations have been discussed in various contexts 
related to quantum groups.  We mention here only a few references:
Cherednik${}^{23}$, Sklyanin${}^{24}$,
Olshanski${}^{6}$, Kulish-Sasaki-Schwiebert${}^{25}$, $\ldots$. 
\par
Suppose that we have an invertible matrix $J\in\EndC(V)$ satisfying 
the reflection equation (3.5). Then we define the matrix 
${\Cal K}\in\EndC(V)\oxC \Uq$ by
$$
{\Cal K}=S({\Cal L^+})JS({\Cal L^-})^t. 
\tag{3.6}
$$
By using this matrix ${\Cal K}=({\Cal K}_{ij})_{ij}$, we define 
the {\it twisted} subalgebra $\Uqk$ of $\Uq$ to be 
the subalgebra generated by the matrix elements ${\Cal K}_{ij}$
($1\le i,j\le N$):
$$
\Uqk=\C[{\Cal K}_{ij} \ \ (1\le i,j\le N)].\tag{3.7}
$$ 
It is easily checked that $\Uqk$ has the coideal property (2.2). 
We remark that the matrix ${\Cal K}$ also satisfies the 
reflection equation similar to (3.5). 
As in Section 2, we denote by ${\Cal J}$ the left ideal associated 
with this subalgebra $\Uqk$:
$$
{\Cal J}=\sum_{1\le i,j\le N} \Uq ({\Cal K}_{ij}-\vep({\Cal K}_{ij})). 
\tag{3.8}
$$
\par
We define a tensor $w_J\in V\oxC V$ of degree 2 associated with $J$
by
$$
w_J=\sum_{1\le i,j\le N}  v_i\ox J_{ij}v_j,\tag{3.9}
$$
where $\{v_j\}_j$ is the canonical basis of $V$.  Then 
the reflection equation (3.5) implies that ${\Cal J}.w_J=0$.  
In this sense, 
the twisted subalgebra $\Uqk$ can be understood as the 
stabilizer of this quadratic tensor $w_J$. 
\par
In the following, we will give a list of constant solutions 
$J$ to (3.5) that we take for the construction of 
$\Uqk$ corresponding to each symmetric space $G/K$. 
$$
\align
\intertext{AI:\qquad\qquad\qquad\quad AII:}
J=&\sum_{k=1}^n e_{kk} a_k
\qquad
J=\sum_{k=1}^n (-e_{2k,2k-1}a_{2k-1}+e_{2k-1,2k}a_{2k})
\quad(a_{2k-1}=q a_{2k}\ \ (1\le k\le n))\\
\intertext{BDI:}
 J=&\sum_{1\leq j,j^\prime \leq \ell} e_{jj} a_j +
 \sum_{\ell<j<\ell^\prime } e_{j j^\prime } q^{-\rho_j }  +
 \sum_{j=1}^\ell e_{jj^\prime} (1-q^{2\rho_\ell} )q^{-\rho_j } 
\ \  (a_1 a_{1^\prime} =\cdots =a_\ell a_{\ell^\prime} =q^{2\rho_\ell})\\  
\intertext{CI:} 
 J=&\sum_{k=1}^{2n} e_{kk} a_k \quad
 (a_1 a_{1^\prime} =\cdots =a_n a_{n^\prime})\\
\intertext{CII:}
J=&\sum_{k=1}^\ell (-e_{2k, 2k-1} a_{2k-1} +e_{2k-1, 2k} a_{2k} 
-e_{(2k-1)^\prime (2k)^\prime} a_{(2k)^\prime } +
   e_{(2k)^\prime (2k-1)^\prime } a_{(2k-1)^\prime } ) \\
   &+\sum_{2\ell<j\leq n } e_{jj^\prime } q^{-\rho_j} 
    -\sum_{2\ell<j\leq n } e_{j^\prime j} q^{\rho_j } 
   +\sum_{j=1}^{2\ell} e_{jj^\prime} (1-q^{2\rho_{2\ell} -2} )q^{-\rho_j } \\
&(a_{2k-1}=qa_{2k} ,~a_{(2k)^\prime} =qa_{(2k-1)^\prime } 
(1\leq k\leq \ell),\ \ a_1 a_{1^\prime} =\cdots =a_{2\ell} 
a_{(2\ell)^\prime} =
 -q^{2\rho_{2\ell} -2}) \\
\intertext{DIII ($n=2\ell$):}
 J=&\sum_{k=1}^\ell (-e_{2k, 2k-1} a_{2k-1} +e_{2k-1, 2k} a_{2k}
-e_{(2k-1)^\prime (2k)^\prime} a_{(2k)^\prime } +
   e_{(2k)^\prime (2k-1)^\prime } a_{(2k-1)^\prime } )\\
  &( a_{2k-1} =q a_{2k}, a_{(2k)^\prime } =qa_{(2k-1)^\prime }
 (1\leq k\leq \ell),\ \  a_1 a_{1^\prime} =\cdots =a_n a_{n^\prime}) \\
\intertext{DIII ($n=2\ell+1$):}
 J =& \sum_{k=1}^\ell (-e_{2k, 2k-1} a_{2k-1} +e_{2k-1, 2k} a_{2k}
-e_{(2k-1)^\prime (2k)^\prime} a_{(2k)^\prime } +
   e_{(2k)^\prime (2k-1)^\prime } a_{(2k-1)^\prime } ) \\
& -e_{n^\prime n} a_n +e_{nn^\prime } a_{n^\prime} \\
&(a_{2k-1} =q a_{2k} ,~a_{(2k)^\prime } =qa_{(2k-1)^\prime }
 (1\leq k\leq \ell), \ \ a_1 a_{1^\prime} =\cdots =a_n a_{n^\prime},
 a_n =a_{n^\prime })
\endalign 
$$
In the list above, the $a_k$'s are nonzero real parameters. 
\Theorem{1}{
\newline{\rm (1) }
Each matrix $J\in\EndC(V)$ listed above satisfies 
the reflection equation (3.5) for the corresponding 
$R$-matrix. 
\newline{\rm (2) }
For any $\ld\in P^+_G$, one has $\dim_\C V(\ld)_{\Cal J}\le 1$. 
Furthermore, the set $P^+_{G,{\Cal J}}$ consisting of all 
$\ld\in P^+_G$ such that $V(\ld)^{\Cal J}\ne 0$ coincides 
with that of the case of $G/K$.  
}
Since $\Uqk$ has the coideal property, we obtain the following 
quantum analogue of the algebra of right $K$-invariant regular 
functions on $G$: 
$$
A_q(G/K):=\{\varphi\in\Aq\,|\, {\Cal J}.\varphi=0\}. \tag{3.10}
$$
Furthermore, as a right $\Uq$-module, this algebra 
has the multiplicity free irreducible decomposition
$$
A_q(G/K)\simeq\bigoplus_{\ld\in P^+_{G,{\Cal J}}} V(\ld)^\vee,
\tag{3.11}
$$
exactly in the same way as in the setting of $G/K$. 
In each case, one can specify an appropriate set of 
parameters $a_k$ so that ${\Cal J}$ associated with $\Uqk$ 
satisfies the condition $S({\Cal J})^\ast\subset {\Cal J}$. 
For such special values of $a_k$, the subalgebra 
$A_q(G/K)$ actually becomes a $\ast$-subalgebra of 
$\Aq$. 
\Remark{1}{
As far as invariant elements are concerned, 
the left ideal ${\Cal J}$ is essential rather than the algebra 
$\Uqk$ itself.
Actually, one can take some different matrices, say
${\Cal K}'$, to get the same left ideal ${\Cal J}$. 
One can take 
${\Cal K}'={\Cal L}^- J^t ({\Cal L^+})^t$,
for instance; for other choices, see Section 2.4 of Noumi${}^8$. 
In our previous paper ${}^8$, we used an {\it additive} way to 
quantize $\frak{k}\subset U(\frak{g})$.  Define 
a matrix ${\Cal M}=({\Cal M}_{ij})_{ij}\in\EndC(V)\oxC\Uq$ by
$$
{\Cal M}={\Cal L}^+- J S({\Cal L}^-)^t J^{-1},\tag{3.12}
$$
and take the coideal
$$
\frak{k}_q=\sum_{1\le i,j\le N} \C {\Cal M}_{ij} \subset\Uq. \tag{3.13}
$$
Then it is easily seen that $\Uq\frak{k}_q={\Cal J}$. 
In this sense, the additive and the multiplicative approaches 
make no essential differences in defining the invariant subalgebra 
$A_q(G/K)$.
}
\Remark{2}{
The definition $\Uq$ of Drinfeld-Jimbo depends on a fixed 
Cartan subalgebra $\frak{t}$ of $\frak{g}$.  
For this reason, we have various possible 
quantizations of the symmetric pair $(\frak{g},\frak{k})$, 
depending on the ``position" of $\frak{k}$ inside $\frak{g}$, 
in relation to the fixed $\frak{t}$. 
We remark that, in the quantum analogues of this section, 
the Lie subalgebra ${\frak{k}}\subset{\frak{g}}$ is taken 
so that,
in the Cartan decomposition $\frak{g}=\frak{k}\oplus\frak{p}$,
$\frak{p}\cap\frak{t}=\frak{a}$ gives a maximal abelian subalgebra 
in $\frak{p}$. 
This choice of $\frak{k}$ is quite opposite 
to the case of $\frak{k}\subset\frak{g}$ corresponding to 
the embedding of root systems.
}
\section{4. Zonal spherical functions}
In the rest of this paper, 
we consider the quantum analogue 
${\Cal H}=A_q(K\backslash G/K)$ 
of the algebra of 
$K$-biinvariant regular functions on $G$, and discuss 
zonal spherical functions associated with finite dimensional
representations of $\Uq$. 
Until now we have checked that the zonal spherical functions 
for the quantum analogue of the following symmetric spaces 
are expressed by Macdonald polynomials associated with 
root systems ${}^{26}$or by Koornwinder's 
Askey-Wilson polynomials for $BC_\ell$ ${}^{27}$:
$$
\alignat 4
\text{(1)} &\quad SU(n)/SO(n) 
 &&\quad  \cdots \quad\text{AI} &&\quad A_{n-1} &&\quad 1 \\
\text{(2)} &\quad SU(2n)/Sp(2n) 
 &&\quad  \cdots \quad\text{AII} &&\quad A_{n-1} &&\quad 4 \\
\text{(3)} &\quad SO(2n)/SO(n)\times SO(n) \ \ (n=\ell)
 &&\quad \cdots \quad\text{DI} &&\quad D_\ell &&\quad 1 \\
\text{(4)} &\quad Sp(2n)/U(n) \ \ (n=\ell)
 &&\quad \cdots \quad\text{CI} &&\quad C_\ell  &&\quad 1,1\tag{4.1}\\
\text{(5)} &\quad Sp(2n)/Sp(2\ell)\times Sp(2\ell)\ \ (n=2\ell)
 &&\quad \cdots \quad\text{CII} &&\quad C_\ell &&\quad 4,3\\ 
\text{(6)} &\quad SO(2n)/U(n) \ \ (n=2\ell) 
 && \quad \cdots \quad\text{DIII} &&\quad C_\ell &&\quad 4,1\\
\text{(7)} &\quad SO(2n)/U(n) \ \ (n=2\ell+1) 
 && \quad \cdots \quad\text{DIII} &&\quad BC_\ell &&\quad 4,1,4
\endalignat
$$
The last two columns of this table indicate the type of the 
restricted root system of each symmetric space and the multiplicity 
of restricted roots $\alpha$ with $(\alpha,\alpha)=2,4,1$.
\par
For the symmetric spaces listed in (4.1), we take the following 
constant solution of the reflection equation:
$$
\alignat 2
J&=\operatorname{diag}(q^{\rho_1},\cdots,q^{\rho_N})
\quad &&\text{for} \quad(1),(3),(4),\\
J&=J_0\operatorname{diag}(q^{\rho_1},\cdots,q^{\rho_N})
\quad &&\text{for} \quad(2),(5),(6),(7),\tag{4.2}\\
\endalignat
$$
where $J_0=\sum_{k=1}^{n}(-e_{2k,2k-1}+e_{2k,2k-1})$. 
The condition $S({\Cal J})^\ast\subset{\Cal J}$ for the left ideal 
${\Cal J}\subset\Uq$ is also fulfilled for this $J$. 
Hence we get the $\ast$-subalgebra
$$
{\Cal H}=A_q(K\backslash G/K)=\{\varphi\in\Aq\,|\, 
{\Cal J}.\varphi=\varphi.\overline{\Cal J}=0 \}.\tag{4.3}
$$
Furthermore, it has the simultaneous eigenspace decomposition 
$$
{\Cal H}=\bigoplus_{\ld\in P^+_{G,{\Cal J}}} {\Cal H}(\ld),
\quad
{\Cal H}(\ld)=\{\varphi\in{\Cal H}\,|\,
C.\varphi=\chi_\ld(C)\varphi\ \ (C\in{\Cal Z}\Uq ) \}
 \tag{4.4}
$$
under the action of the center ${\Cal Z}\Uq$,
where $\dim_\C {\Cal H}(\ld)=1$ for all $\ld\in P^+_{G,{\Cal J}}$.
\par
As for the cases (1) and (2), a detailed description of the 
structure of ${\Cal H}$ and the zonal spherical functions 
is already given in Ueno-Takebayashi${}^{28}$ and Noumi ${}^8$.  
Hereafter, we will consider 
the remaining cases (3) -- (7) , so that $G=SO(2n)$ or $G=Sp(2n)$.
In order to describe the algebra ${\Cal H}$, recall that 
the quantum group $G_q$ has the diagonal subgroup ${\Bbb T}$ which
is isomorphic to the $n$-dimensional torus.  Namely, we have 
a canonical surjective $\ast$-homomorphism 
$\cdot|_{\Bbb T} : \Aq\to A({\Bbb T})$,
where $A({\Bbb T})=\C[z_1^{\pm 1},\cdots,z_n^{\pm 1}]$ with 
$\ast$-operation $z_j^\ast=z_j^{-1}$.
By means of this ``restriction mapping", we consider the 
composition 
$$
\cdot|_{\Bbb T} : 
 {\Cal H}=A_q(K\backslash G/K) \hookrightarrow 
\Aq\twoheadrightarrow A({\Bbb T}).\tag{4.5}
$$
Let $W(\Sigma)$ be the Weyl group of the restricted root 
system $\Sigma$ of the symmetric space $G/K$, 
and $P(\Sigma)$ the lattice of integral weights 
of $\Sigma$. 
We define the elements $x_1,\cdots,x_\ell$ in 
$A({\Bbb T})=\C[z_1^{\pm 1},\cdots,z_n^{\pm 1}]$ 
by
$$
\alignat 2
&x_j=z_j^2\quad(1\le j\le\ell)&&\quad\text{for}\quad (3),(4),\\
&x_j=z_{2j-1}z_{2j}\quad(1\le j\le\ell)&&\quad\text{for}\quad (5),(6),(7),
\tag{4.6}
\endalignat
$$
and take the following subalgebra ${\Cal A}$ of 
$A({\Bbb T})=\C[z_1^{\pm 1},\cdots,z_n^{\pm 1}]$:
$$
{\Cal A}=\bigoplus_{\mu\in P(\Sigma)}\,\C\,x^\mu.  \tag{4.7}
$$
(We understand $x_j^\half=z_j$ in the case (3).)
\Theorem{2}{
\newline{\rm (1)}
The restriction mapping $\cdot|_{\Bbb T} : {\Cal H}\to A({\Bbb T})$
is injective.  Hence ${\Cal H}\subset \Aq$ is a commutative 
subalgebra.
\newline{\rm (2)}
The image ${\Cal H}|_{\Bbb T}$ of the restriction mapping 
coincides with the subalgebra ${\Cal A}^{W(\Sigma)}$ 
consisting of $W(\Sigma)$-invariant elements in ${\Cal A}$ above. 
Namely we have ${\Cal H}\simeq {\Cal A}^{W(\Sigma)}$.
}
We recall now the Macdonald polynomials, restricting ourselves 
to those associated with root systems of type $B_\ell,C_\ell,D_\ell$. 
Taking a root system $\Sigma$ of type $B_\ell,C_\ell,D_\ell$,
let $(m_\alpha)_{\alpha\in\Sigma}$ be a set of nonnegative real 
numbers, invariant under the action of the Weyl group 
$W=W(\Sigma)$ of $\Sigma$. 
Setting $t_\alpha=q_{\alpha}^{m_\alpha/2}$ with 
$q_\alpha=q^{(\alpha,\alpha)/2}$, 
we take the function $\Dt^+(x)$ on $(\C^\ast)^\ell$ defined by
$$
\Dt^+(x)=\prod_{\alpha\in\Sigma^+}
\frac{(x^\alpha;q_{\alpha})_{\infty}}
{(t_{\alpha}x^{\alpha};q_{\alpha})_{\infty}},\tag{4.8}
$$
where $(a;q)_\infty=\prod_{i=0}^\infty(1-aq^i)$, 
and set $\Dt(x)=\Dt^+(x)\,\Dt^+(x)^\ast$.
With the $\ast$-operation $x_j^\ast=x_j^{-1}$, define the 
scalar product $\langle{\,\,|\,\,}\rangle$ on 
the subalgebra 
${\Cal A}^{W(\Sigma)}$ of
$W(\Sigma)$-invariants in 
${\Cal A}=\bigoplus_{\mu\in P(\Sigma)} \C x^\mu$
by
$$
\langle{\,f|g\,}\rangle=\frac{1}{|W(\Sigma)|}
\int_T f(x)^\ast\,g(x)\Dt(x)
\quad (f,g\in {\Cal A}^{W(\Sigma)}), 
\tag{4.9}
$$
where $\int_T$ is the normalized Haar measure of the 
torus 
$T={\Bbb R}^\ell/Q^\vee(\Sigma)$, and 
$x^\mu$ are regarded as functions on $T$ by 
$x^\mu(\nu)=e^{2\pi\sqrt{-1}(\mu,\nu)}$. 
Then it is known by Macdonald${}^{26}$  that 
the subalgebra ${\Cal A}^{W(\Sigma)}$ 
has a unique basis 
$\{P_\mu(x)\}_{\mu\in P^+(\Sigma)}$, orthogonal under the 
scalar product (4.9), such that
$$
P_\mu(x)=m_\mu+{\sum}_{\nu < \mu } u_{\mu\nu}m_\nu(x) \tag{4.10}
$$
for each $\mu\in P^+(\Sigma)$, where 
$m_\mu(x)=\sum_{\nu\in W(\Sigma)\mu} x^\nu$ is the orbit sum and 
$<$ denotes the dominance order of weights.
These $P_\mu(x)$ ($\mu\in P^+(\Sigma)$) are called 
the {\it Macdonald polynomials} associated with the root 
system $\Sigma$ (or, to be more precise, associated with the pair of root systems $(\Sigma,\Sigma^\vee)$). 
An extension of Macdonald polynomials of type $BC_\ell$ has 
been introduced by Koornwinder ${}^{27}$.  
Koornwinder's {\it Askey-Wilson polynomials for $BC_\ell$} 
are defined similarly by taking the root system $\Sigma$ of 
type $BC_\ell$ and the function 
$$
\Dt^+(x)=\prod_{k=1}^\ell 
\frac{(x_k^2;q)_\infty}{(ax_k,bx_k,cx_k,dx_k;q)_\infty}
\prod_{1\le i<j\le \ell} 
\frac{(x_i/x_j,x_ix_j;q)_\infty}
{(tx_i/x_j,tx_ix_j;q)_\infty},\tag{4.11}
$$
instead of $\Dt^+(x)$ of (4.8), where $a,b,c,d,t$ are real parameters
and $(a,b,\cdots;q)_\infty=(a;q)_\infty (b;q)_\infty\cdots$.
Among many works related to Macdonald polynomials, we only refer to 
Cherednik${}^{29}$
for the relation to affine Hecke algebras and to 
Etingof-Kirillov${}^{30}$
for a realization of Macdonald polynomials of type $A$ by 
vector-valued characters of $U_q(\frak{gl}(N))$. 
\par\medskip
We now return to the setting of quantum symmetric spaces
of type (4.1), (3)--(7) and 
consider the zonal spherical functions $\varphi_\ld$ 
($\ld\in P^+_{G,{\Cal J}}$) on 
$(G/K)_q$ defined by the constant solution $J$ of (4.2). 
{}From our definition of $\Uqk$, 
it turns out that the zonal spherical 
function $\varphi_{\ld}$ 
can be normalized so that its restriction to ${\Bbb T}$ 
takes the form 
$$
\varphi_{\ld}|_{\Bbb T}=z^{\ld}+ {\sum}_{\mu<\ld} a_{\ld\mu}z^{\mu},
\tag{4.12}
$$
where $<$ stands for the dominance order of weights.
One can show that the restriction $\varphi_{\ld}|_{\Bbb T}$ 
are expressed by Macdonald polynomials (or Koornwinder's
Askey-Wilson polynomials for $BC_\ell$) in the variables 
$x=(x_1,\cdots,x_\ell)$ of (4.6). 
For $\ld\in P^+_{G,{\Cal J}}$, we take the dominant 
integral weight $\mu\in P^+(\Sigma)$ of the corresponding 
restricted root system such that  $z^\ld=x^\mu$ under (4.6).
Then we have
\Theorem{3}{ 
If the symmetric space $G/K$ is of type (3)--(6), 
then the restriction $\varphi_\ld|_{\Bbb T}$ 
($\ld\in P^+_{G,{\Cal J}}$ ) of the zonal spherical 
function on our quantum symmetric space 
$(G/K)_q$ coincides with the Macdonald 
polynomial $P_\mu(x)$ associated 
with the restricted root system $\Sigma$ and 
and the set $(m_\alpha)_{\alpha\in\Sigma}$ 
of multiplicities of roots of $G/K$.  
They have the base $q^4$ in the cases (3),(4) and $q^2$ 
in the cases (5),(6), respectively. 
If $G/K$ is of type (7), $\varphi_\ld|_{\Bbb T}$ coincides 
with Koornwinder's Askey-Wilson polynomial $P_\mu(x)$ 
for $BC_\ell$ with parameters $(a,b,c,d; q,t)$ replaced by 
$(q^3,q^3,-q,-q;q^2,q^4)$.
}
\par
By means of the function $\Dt^+(x)$ of (4.8) (or (4.11)), 
we define the $q$-difference operator $D_\sigma$ 
for the weight $\sigma=\epsilon_1$ as follows:
$$
D_\sigma=\frac{1}{|W_\sigma|}
\sum_{w\in W} w\Phi_\sigma(x)(T_{w\sigma}-1)\quad\text{with}\quad
\Phi_\sigma(x)=\frac{T_\sigma\Dt^+(x)}{\Dt^+(x)},\tag{4.13}
$$
where $T_\mu$ is the $q$-shift opertar such that 
$T_\mu x^{\ld}=q^{(\ld,\mu)}x^\mu$.  
This $q$-difference operator $D_\sigma$ is selfadjoint 
with respect to the scalar product of Macdonald (or Koornwinder). 
Furthermore, the Macdonald polynomials (or Koornwinder's 
Askey-Wilson polynomials) $P_\mu(x)$ are  eigenfunctions 
of the operator $D_\sigma$.  
\par
In our context, 
the operator $D_\sigma$, with the parameters described in Theorem 3, 
arises as the radial component 
of the central element 
$$
C_\sigma=\sum_{1\le i,j\le N} q^{2\rho_i}L^+_{ij}S(L^-_{ji})\tag{4.14}
$$
of $\Uq$ defined by Reshetikhin-Takhtajan-Faddeev ${}^2$. 
For a central element $C\in{\Cal Z}\Uq$, we mean by 
the {\it radial component} of $C$ an operator 
$D:{\Cal A}^{W(\Sigma)}\to {\Cal A}^{W(\Sigma)}$ that makes the 
following diagram commutative:
$$
\CD
{\Cal H} @>C.>> {\Cal H}\\
@V{\cdot|_{\Bbb T}}VV  @VV{\cdot|_{\Bbb T}}V \\
{\Cal A}^{W(\Sigma)} @>>D> {\Cal A}^{W(\Sigma)}.
\endCD
\tag{4.15}
$$
In fact it turns out that 
the radial component of $C_\sigma-\chi_{\ld}(C_\sigma)$ 
is a constant multiple of an operator of the 
form $D_\sigma - a_\sigma(\mu)$  for some $a_\sigma(\mu)\in \C$. 
Hence the restriction $\varphi_\ld |_{\Bbb T}$ is an 
eigenfunction of $D_\sigma$ for each $\ld\in P^+_{G,{\Cal J}}$.  
By using this fact, one can show that, on the $\ast$-subalgebra 
${\Cal H}=A_q(K\backslash G/K)$, the scalar product defined 
by the Haar functional of $G_q$ in (1.6) corresponds to 
that of Macdonald (or Koornwinder) in (4.9), up to a scalar multiple. 
{}From this it follows that the basis 
$\{\varphi_\ld|_{\Bbb T}\}_{\ld\in P_{G,{\Cal J}}}$ 
of ${\Cal A}^{W(\Sigma)}$ coincides
with the orthogonal basis of Macdonald (or Koornwinder) described 
above. 
%

%
%
\par\bigskip\par\noindent
{\bf Notes:}
This article is based on the talk given at the conference by 
one of the authors (M.N.).  On that occasion, Prof. Allan Solomon 
kindly made the following two comments: 
\roster
\item It is a fairly strong abuse of language to describe as a 
``quantum group" an object which is not a group; but you 
have defined as ``quantum subgroups" objects which are not 
even quantum groups!
\item The commutation $ef-fe=(t-t^{-1})/(q-q^{-1})$ would 
indicate a Haar measure based on the $q$-integral of 
Nelson et al. (1988) and not that of Jackson (1909). 
\endroster
\enddocument